# Tunable band gaps of axially moving belt on periodic elastic foundation

Lei Lu[*]

*Department of Mathematics, Shanghai Lixin University of Accounting and Finance, Shanghai 201209, China*

**Abstract**

The present paper investigates the band structure of an axially moving belt resting on a foundation with periodically varying stiffness. It is concluded that the band gaps appear when the divergence of the eigenvalue occurs and the veering phenomenon of mode shape begins. The bifurcation of eigenvalues and mode shape veering lead to wave attenuation. Hence, the boundary stiffness modulation can be designed to manipulate the band gap where the vibration is suppressed. The contribution of the system parameters to the band gaps has been obtained by applying the method of varying amplitudes. By tuning the stiffness, the desired band gap can be obtained and the vibration for specific parameters can be suppressed. The current study provides a technique to avoid vibration transmission of the axially moving material by designing the foundation stiffness.

## 1. Introduction

The axially moving belts are widely employed in various engineering applications, such as driving belts, transmission belts, magnetic belts, et al. The periodic structure is also a common device that can be found in modern industrial fields as vibration isolator and mechanical filter. The dynamic properties of periodic structures, similar to phononic crystals, are attractive to engineers and researchers. The question then arises:

---

[*] E-mail address: luleic@163.com

Is that possible the vibration of axially moving materials can be suppressed by the periodic designed foundation instead of constant foundation?

The axially moving continua have been studied for centuries [1]. The transverse vibrations have been investigated widely [2-5]. The analysis and control of transverse vibration of axially moving string was studied systematically by Chen [6]. Considering various constraints [7-9], the axially moving continua were investigated in the available literature. Yang et al. [10] investigated the axially moving beam resting on elastic foundation with constant stiffness and obtained the explicit expression of critical velocity. Bhat et al. [11] investigated the dynamic response of the moving belt supported by elastic foundation by numerical techniques. The researches of the foregoing literature are about the uniform structures. When the systems involve periodic parameters, such as cross-section [12], stiffness, etc., the systems transform into periodic structures. One of the distinctive features of the periodic structures is the appearance of band gaps. The wave attenuation occurs in the band gaps for systems composed of periodic units [13]. The fascinating characteristics of the band gaps of the periodic structures have drawn much attention for the recent decades.

The properties of the wave propagation in periodic structures were investigated by Brillouin [14] systematically. Mead [15] reviewed a notion of the wave propagation constants in periodically supported infinite beams. Asfar and Nayfeh [16] applied the method of multiple scales to the periodic structures. The axially moving belt with periodic cross-section has been investigated by Sorokin [17] recently. He applied the method of varying amplitude to the periodic structures for the vibration suppression. Inspired by Sorokin [18], a uniform axially moving belt with periodic varying foundation stiffness and stiffness modulation which composes a periodic structure is considered. The band gap depends on the system parameters which would vary in different cases. The investigation leads to a tunable band gap that stems from the research of photonic crystal [19]. Some investigations of the tunable band gap focused on the photonic crystal or phononic crystal [20-22]. Various wave guides were designed to control the wave propagation or tune the band gap [23-26]. Lv et al. [27] demonstrated a continuous tunable band gap of a metamaterial beam by modulating the moment of inertia of the beam. In most recently, Yang et al. [28] investigated a two-dimensional periodic lattice tuned by parametric excitations to control the wave propagation.

The wave properties of an axially moving belt resting on elastic foundation with periodic stiffness are examined in the paper. It is found that the vibration can be suppressed by designing specific periodic foundation stiffness. The rest of the work is structured as follows. **Section 2** establishes the governing equation for the system with periodic structure. **Section 3** calculates the solution of the governing equation by method of varying amplitude, and obtains the analytical solution. **Section 4** presents the dispersion relation and the band gaps of the model by employing the eigenfunction, and interprets the mechanism of the wave attenuation in term of the mode shapes aberration. **Section 5** describes the tunable band gap varying with the parameters of the system and reveals the relation of band gap with the foundation stiffness, the stiffness modulation and the axial velocity. **Section 6** ends the work by drawing the main conclusions.

## 2. Governing equations

Consider an axially moving belt with velocity $V$ resting on a periodic elastic foundation as shown in **Fig. 1**. The displacement in the transverse direction is denoted by $U(X, T)$, where $X$ is the spatial coordinate and $T$ is time. The linear density is $\rho$ and the tension is $P$. $S(X)$ denotes the varying stiffness of the elastic foundation with the period $\Phi$ at the point $X$, and $S(X)= S(X+\Phi)$.

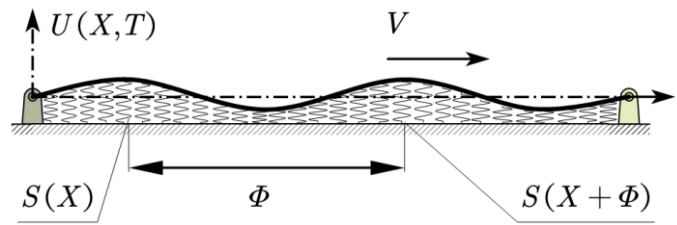

Fig. 1 The schematic of axially moving belt resting on a continuous periodic elastic foundation.

The governing equation is derived instantly by Newton's second law or Hamilton's

principle as

$$\rho(U_{TT} + 2VU_{XT} + V^2 U_{XX}) - PU_{XX} + S(X)U = 0. \tag{1}$$

Assuming that the stiffness $S(X)$ varies harmoniously

$$S(X) = S_0(1 + \sigma \cos \varphi X), \tag{2}$$

where $S_0$ is a constant that stands for the stiffness of foundation without modulation, $\sigma$ represents the stiffness modulation, $\varphi = 2\pi/\Phi$ denotes the density of the stiffness varying with length which represents one-dimensional reciprocal lattice [14].

Introducing the nondimensional variables and parameters

$$u = \frac{U}{\Phi}, \varphi = \frac{2\pi}{\Phi}, c = \sqrt{\frac{P}{\rho}}, x = \varphi X, t = \varphi c T, v = \frac{V}{c}, s = \frac{S_0}{\varphi^2 P} \tag{3}$$

Eq. (1) is rewritten as

$$u_{tt} + 2v u_{xt} - (1 - v^2) u_{xx} + s(1 + \sigma \cos x) u = 0. \tag{4}$$

The solution to Eq.(4) is assumed as

$$u(x,t) = A(x) e^{i\omega t}. \tag{5}$$

By substituting Eq. (5) into Eq. (4), one obtains

$$(1 - v^2) A'' - 2i\omega v A' + [\omega^2 - s(1 + \sigma \cos x)] A = 0. \tag{6}$$

When the foundation stiffness $S(x)=0$, the system degenerates into the same form in [2].

## 3. Solution by the method of varying amplitudes

By applying the method of varying amplitudes [17, 18], the author seeks the solution to Eq. (6) in the form of an infinite series as

$$A(x) = \sum_{m=-\infty}^{+\infty} a_m(x) e^{imx}. \tag{7}$$

Substituting Eq. (7) back into Eq. (6), one obtains

$$\sum_{m=-\infty}^{+\infty} \left\{ (1-v^2)(a_m'' + 2ima_m' - m^2 a_m) - 2i\omega v(a_m' + ima_m) + \left[\omega^2 - s(1+\sigma \cos x)\right] a_m \right\} e^{imx} = 0$$

$$. \tag{8}$$

By transforming cos*x* to the well-known Euler's formula and balancing the coefficients of the harmonic involved, a set of equations is obtained as

$$(1-v^2)(a_m'' + 2ima_m' - m^2 a_m) - 2i\omega v(a_m' + ima_m) + (\omega^2 - s)a_m - \frac{1}{2}s\sigma(a_{m-1} + a_{m+1}) = 0, m \in \mathbb{Z}. \tag{9}$$

Assuming $K=-ik$ is the eigenvalue of Eq. (9), where $k$ is the wave number, and $\boldsymbol{\alpha} = (\alpha_0\ \alpha_1\ \alpha_{-1}\ \alpha_2\ \alpha_{-2}\ \ldots)^T$ is the associated eigenvector. The coefficients of the right-hand side of Eq. (7) can be rewritten in vector form as

$$\alpha(x) = \begin{pmatrix} a_0(x) \\ a_1(x) \\ a_{-1}(x) \\ a_2(x) \\ a_{-2}(x) \\ \vdots \end{pmatrix} = \begin{pmatrix} \alpha_0 \\ \alpha_1 \\ \alpha_{-1} \\ \alpha_2 \\ \alpha_{-2} \\ \vdots \end{pmatrix} e^{-ikx}. \tag{10}$$

Then Eq. (9) turns to the matrix form

$$\mathbf{A}\alpha = 0, \tag{11}$$

where the matrix **A** has the following elements

$$a_{mn} = \begin{cases} (v(k-m) - \omega)^2 - (k-m)^2 - s, & m = n; \\ -\dfrac{1}{2}s\sigma, & m = n \pm 1; \\ 0, & \text{otherwise.} \end{cases} \tag{12}$$

Here, $a_{mn}$ denotes the element in the $m^{th}$ row and the $n^{th}$ column of the matrix **A**.

Combining Eqs. (5), (7) and (10), the solution to Eq. (4) can be written as

$$u(x,t) = \sum_{m=-\infty}^{+\infty} \alpha_m e^{i((m-k)x + \omega t)}, \tag{13}$$

where the frequency is $\omega$, and the compound wavenumber is $m$-$k$. From Eq. (13) it follows that for a continuous periodic system, the corresponding compound wave comprises infinite components with the same frequency, but different wavenumbers.

## 4. Dispersion relation and band gaps

## 4.1 Dispersion relation

This section is devoted to revealing the mechanism of the dispersion relation and the band gaps by studying the eigenvalues of the system. To obtain the non-trivial roots of Eq. (11), the determinant of the matrix **A** should be equal to zero, which leads to the following eigenfunction

$$\begin{aligned}
f(\omega,k,v,s,\sigma) = det(\mathbf{A}) &= \prod_{m=-\infty}^{+\infty}[(v(k-m)-\omega)^2-(k-m)^2-s] \\
&-(\frac{1}{2}s\sigma)^2 \sum_{n=-\infty}^{+\infty} \prod_{m \neq n, n+1}^{+\infty}[(v(k-m)-\omega)^2-(k-m)^2-s] \\
&+(\frac{1}{2}s\sigma)^4 \sum_{\substack{n=-\infty \\ l \geq n+2}}^{+\infty} \prod_{m \neq n, n+1, l, l+1}^{+\infty}[(v(k-m)-\omega)^2-(k-m)^2-s] \\
&-(\frac{1}{2}s\sigma)^6 \sum_{\substack{i=-\infty \\ n \geq i+2 \\ l \geq n+2}}^{+\infty} \prod_{m \neq i, i+1, n, n+1, l, l+1}^{+\infty}[(v(k-m)-\omega)^2-(k-m)^2-s] + \cdots \\
&= 0,
\end{aligned} \quad (14)$$

where $f$ is an implicit function with five parameters: frequency $\omega$, wavenumber $k$, axial velocity $v$, foundation stiffness $s$ and stiffness modulation $\sigma$. The dispersion relation $\omega(k)$ can be determined by the given parameters. Thus, Eq. (14) can be rewritten as an implicit function with respect to $\omega$ and $k$, i.e., $f(\omega, k)=0$. For any integer $m$, $f(\omega, k+m) = f(\omega, k)$, i.e., $f(\omega, k)=0$ is an implicit periodic function of wavenumber $k$ with the period 1, where 1 is the least positive period of $f$. Uncertainty in the wavenumber $k$ is best avoided by restricting $k$ to the interior of the first Brillouin zone, i.e., $k \in [-0.5, 0.5]$.

In order to investigate the dispersion relation, the frequency and the wavenumber should be drawn out from eigenfunction. It is obvious that Eq. (14) is a function of $(s\sigma)^{2m}$, where $m$ is an integer. If $s<1$ and $\sigma<1$, then $(s\sigma)^{2m} \ll 1$, thus, the result is reliable when the higher order terms of the function are neglected. As an illustration, the author starts on the first band gap. Assuming that the foundation stiffness $s=0$, or $s \ll 1$, it is reasonable to take into account only one term of the eigenfunction, i.e., let $m=0$ in Eq. (14). Then, it follows

$$f = (vk-\omega)^2 - k^2 - s = 0. \quad (15)$$

According to the definition of the group velocity $d\omega/dk$ [29], one obtains

$$c_g = \frac{d\omega}{dk} = v \pm \frac{k}{\sqrt{k^2 + s}}. \tag{16}$$

The direction of wave propagation can be detected by the sign of the group velocity, namely, the positive sign denotes forward travelling wave and the negative one denotes backward travelling wave. Thus, with regards to the axially moving belt, the corresponding compound waves are divided into two sets with different wavenumbers but the same frequency. However, each frequency corresponds to one wavenumber in the first Brillouin zone.

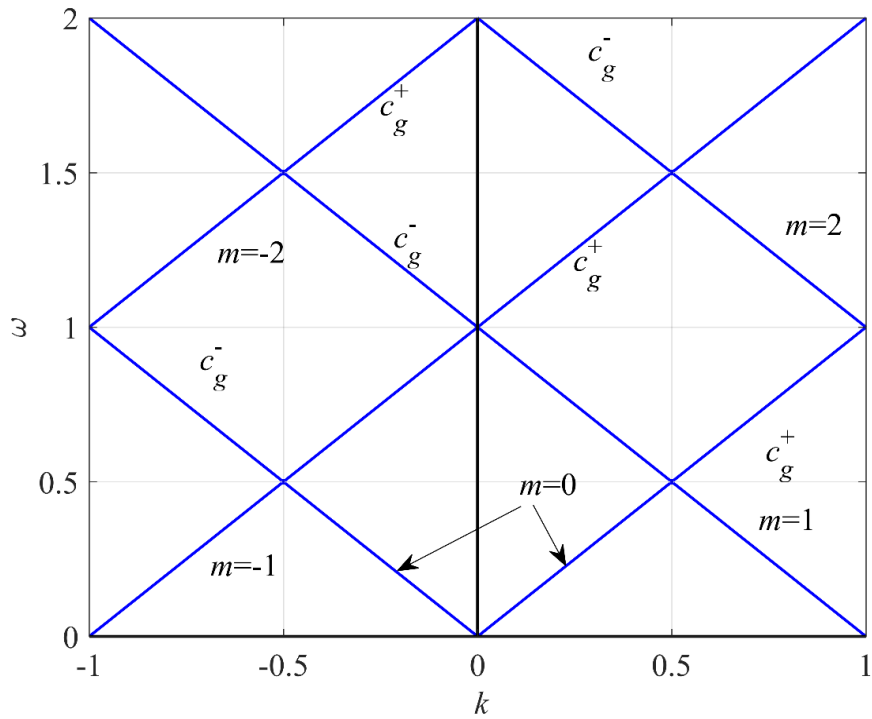

(a)

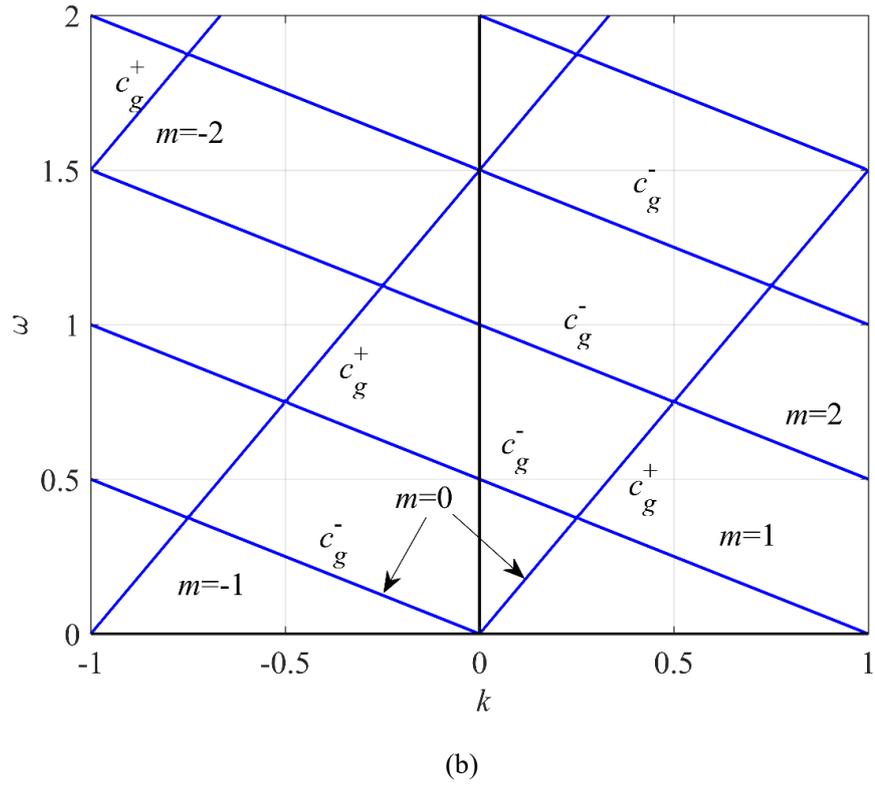

(b)

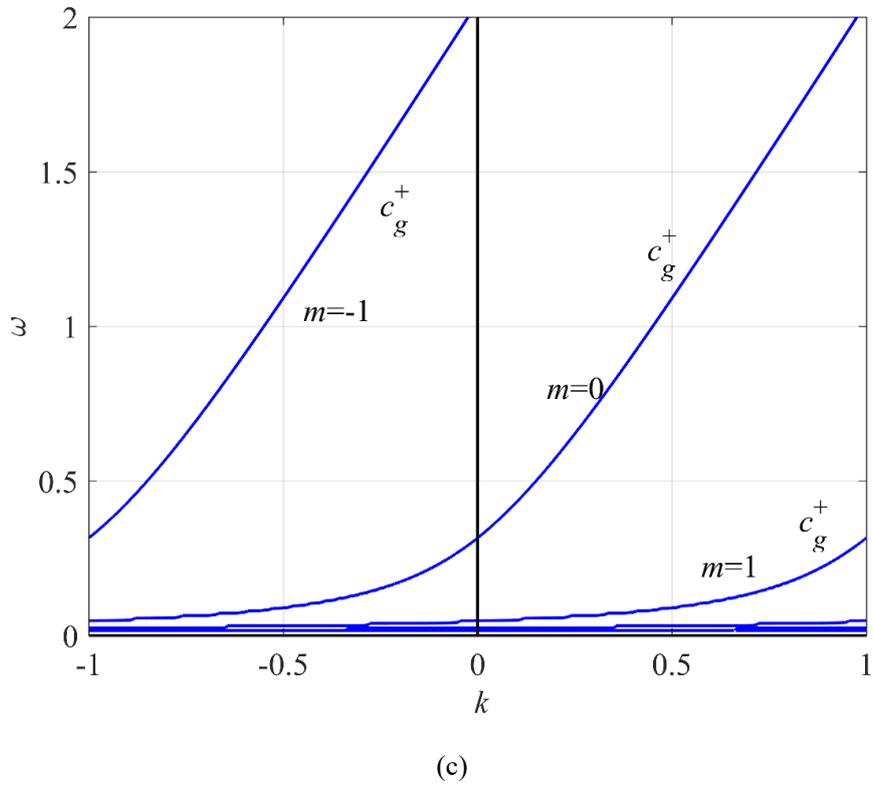

(c)

**Fig. 2** Dispersion relation $\omega(k)$, for (a) $v=0$, $s=0$; (b) $v=0.5$, $s=0$ and (c) $s=0.1$, $\sigma=0.5$, $v=0.999$ (the critical velocity is 1 [2]).

As an illustration, **Fig. 2** demonstrates the dispersion relation of the belt with different parameters in the case of band gaps vanishing. If the foundation stiffness decreases toward zero, the band gap vanishes, which agrees well with the result of [17] as shown in **Fig. 2** (a). Where $m$ represents the $m^{th}$ order mode, $c_g^+$ and $c_g^-$ denote the group velocities of forward travelling wave and backward wave, respectively. The compound wave of the system comprises with different modes which can be divided into two opposite directions. The slope of the lines denotes the group velocity of the wave propagation. The positive one represents the forward travelling wave and the negative one represents the backward travelling wave. The speed of the forward travelling wave increases with the axial velocity while the backward travelling wave speed decreases with the axial velocity, as demonstrated in **Fig. 2** (b). When the axial velocity limits to the critical velocity, the slop of the lines, i.e. the speed of the backward travelling wave approaches zero as shown in **Fig. 2** (c). The maximum of speed of the forward travelling wave is $v+1$ which can be derived from Eq. (16) instantly, where $v$ is the axial velocity.

**4.2 The existence of band gaps**

The dependence of eigenvalues on system parameters can be illustrated by a family of loci [30]. With the varying of the foundation stiffness, the frequency curves veer and eigenvalues diverge. Analogously, the eigenfunction mode shapes veering in the transition zone [31], the author finds that the mode shapes of the wave change dramatically in the band gap. As shown in **Fig. 3(a)**, when $\omega = 0.28$, $\omega = 0.4$, $\omega = 0.5$, mode shapes in the band pass are normal, when $\omega = 0.27$ and $\omega = 0.45$, mode shapes are irregular in the first band gap in **Fig. 3(b)** and the second band gap in **Fig. 3(c)**, respectively. The mode shapes undergo abrupt changes when the frequency changes from band pass to band gap. If the length of the compound wave modes becomes infinite, the wavenumber or reciprocal wave length approaches zero, i.e., the wave attenuates in these areas. The survey of wave propagation in the vicinity of band gap shows that: the modes of wave are normal when the frequency falls in the band pass, while the modes become abruptly and the bounded amplitude is broken when the frequency falls in the band gap, where the wave propagation is stopped.

Another interpretation of the wave attenuation is in terms of the wave velocity. When the speed of forward travelling wave is equal to the backward one, the superposition of wave velocity is zero, i.e., the wave propagation is stopped. These points in frequency-wavenumber plane are named as *veering points*. The mode shapes veering and the eigenvalues changing occur instantaneously at these points. As demonstrated in **Fig. 4 (a)**, the slopes of the curves denote the group velocities of the wave, and $c_g^+$ and $c_g^-$ represent the same meaning as the previous section.

Actually, the number of the band gaps is dependent on the terms of the expansion in Eq. (14), and the maximum number of the band gaps is $2m+1$. If $m=0$, then $2m+1=1$, i.e. the maximum number of the band gaps is one. If the expansion Eq. (14) is restricted in the first Brillouin zone, let $\sigma=0$, Eq. (14) degenerates to the same case of $m=0$. The dispersion relation is as shown in **Fig. 4 (a)** which has only one band gap. In fact, when $m=1$, there will be three band gaps as shown in the remaining of **Fig. 4**, but it is too tiny to be observed. If $m$ is not restricted, there are also fourth band gap or higher band gaps. Sorokin and Thomsen [17] analyzed the band gap by employing the perturbation method, they found only one band gap because they truncated the perturbation parameter in second order, so they missed the higher frequency band gaps.

Magnifying the third band gap and fourth band gap of **Fig. 4 (b)** a thousand and five thousand times respectively yields the enlarged view as shown in **Fig. 5**. Three different values of $\sigma$ are given, which demonstrates that the stiffness modulation effects on the second or higher band gap. If $\sigma$ decreases toward zero, the width of the second or higher band gap approaches zero. Thus, the existence and the number of the band gap are dependent upon the periodic foundation stiffness and the stiffness modulation of the axially moving belt. The system, at least, has one band gap if the foundation stiffness exists. There will be more than two band gaps if the stiffness modulation appears with the foundation stiffness. The higher band gaps accompany with foundation stiffness and the stiffness modulation.

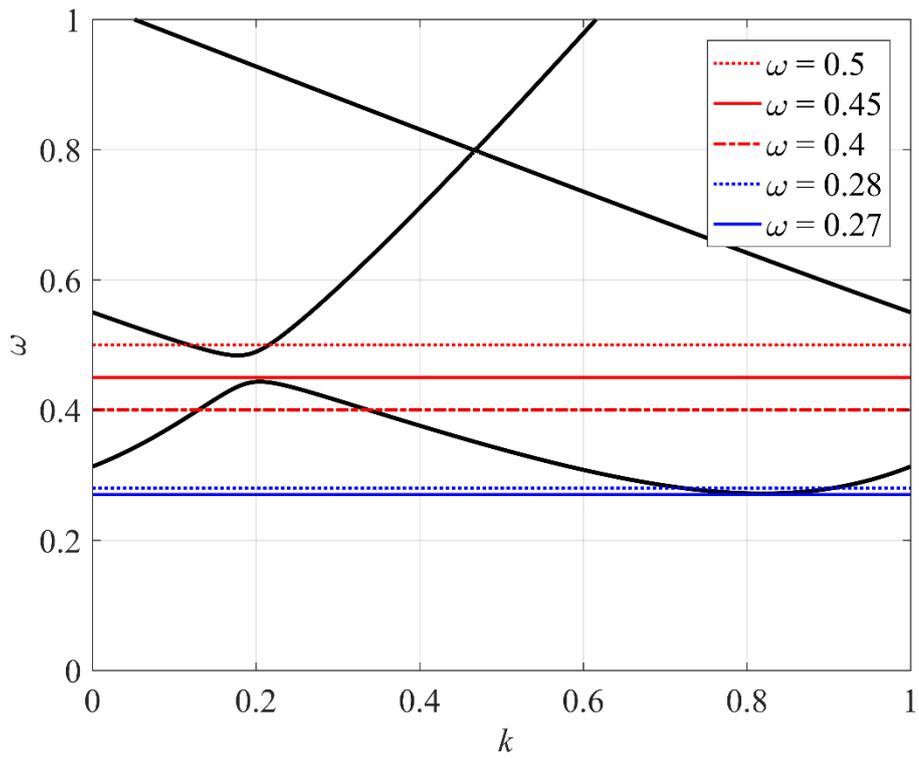

(a)

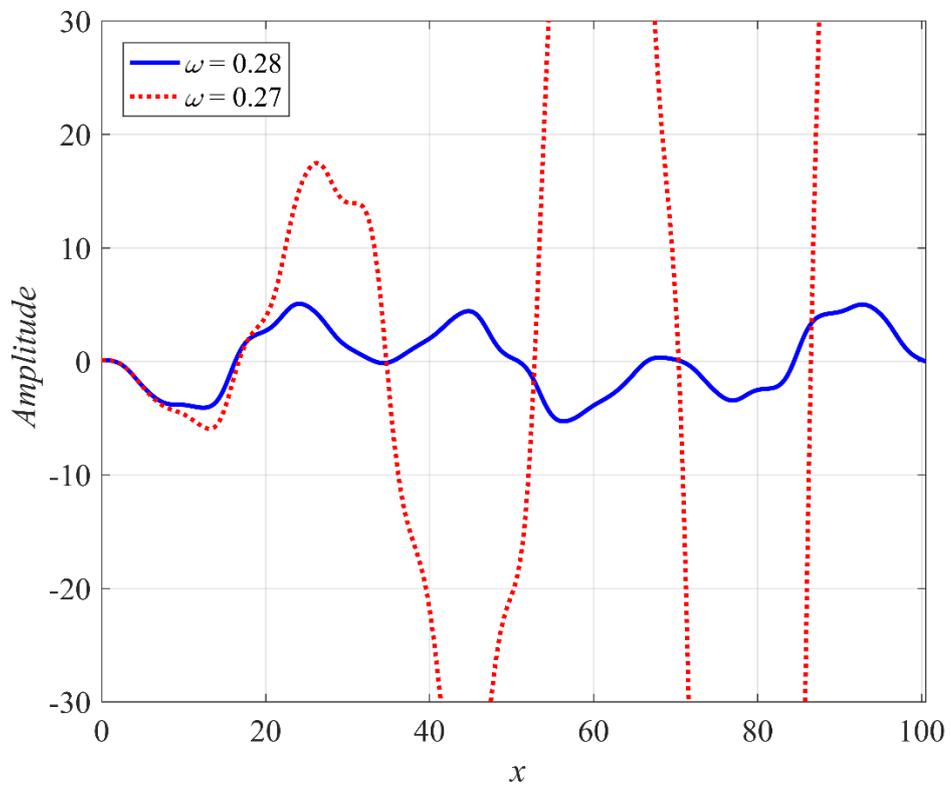

(b)

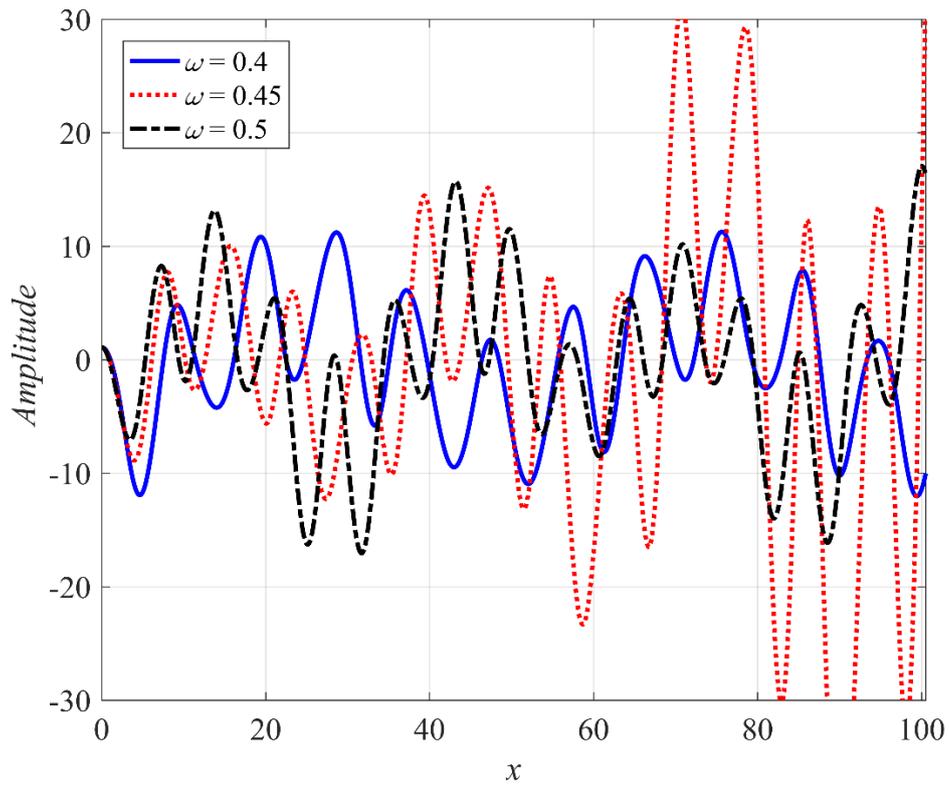

(c)

**Fig. 3**. Mode shapes change in the band gaps, for $v$=0.5, $s$=0.1, $\sigma$=0.5. (a) frequencies in the vicinity of two band gaps; (b) mode shapes change in the first band gap; (c) mode shapes change in the second band gap.

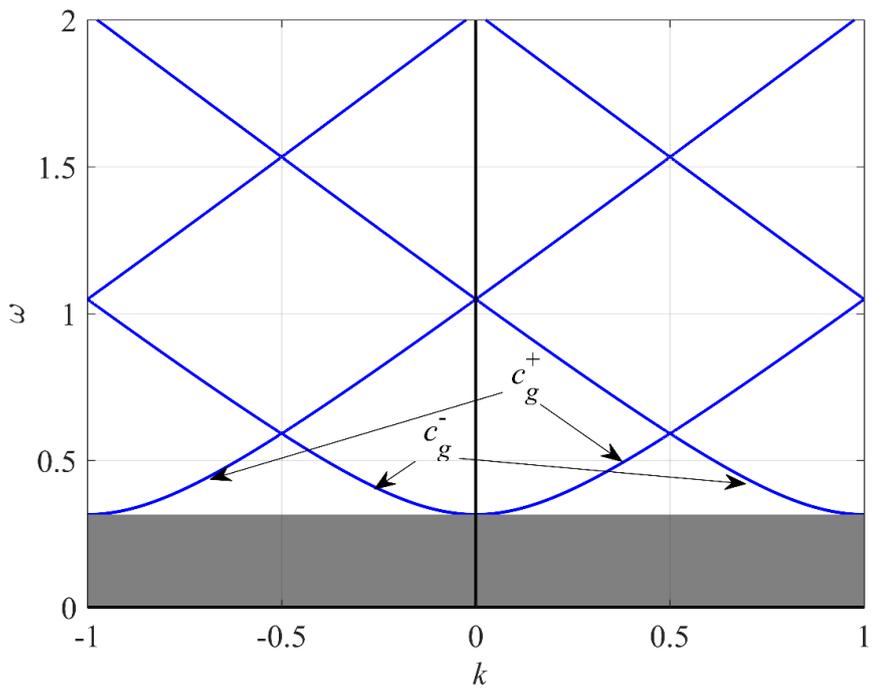

(a)

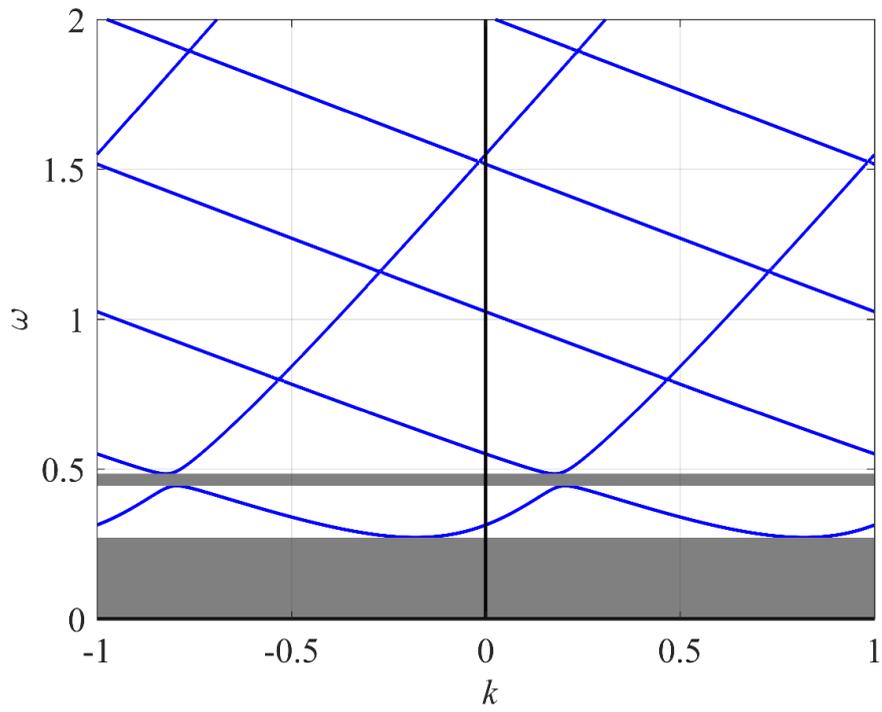

(b)

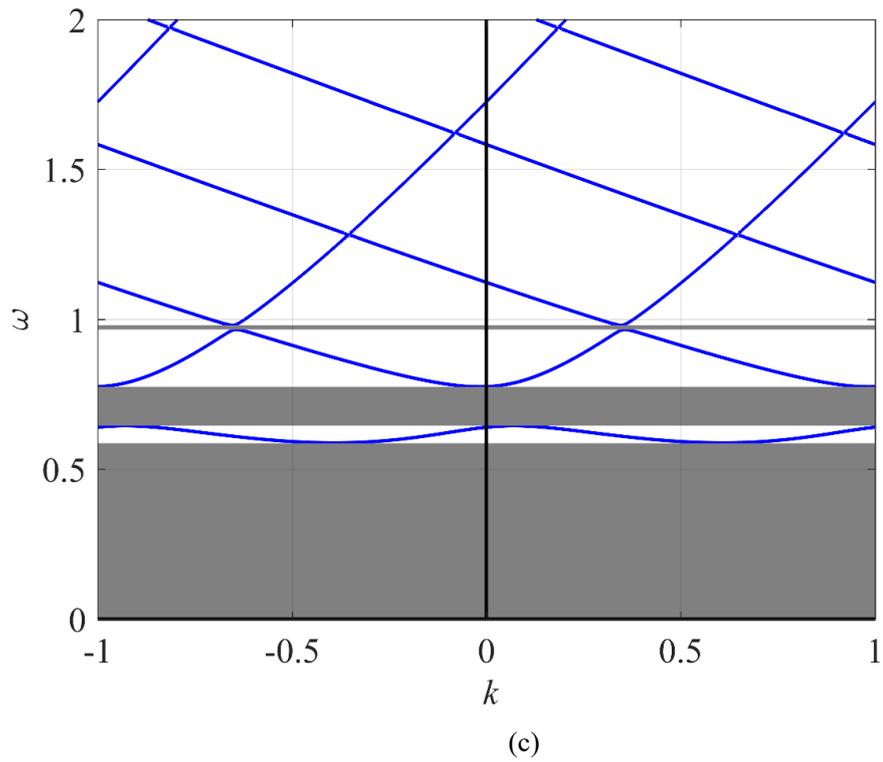

(c)

**Fig. 4** The band gaps are marked by gray shaded areas. (a) $s$=0.1 $v$=0, $\sigma$=0, (b) $s$=0.1 $v$=0.5, $\sigma$=0.5, and (c) $s$=0.5, $v$=0.5, $\sigma$=0.5.

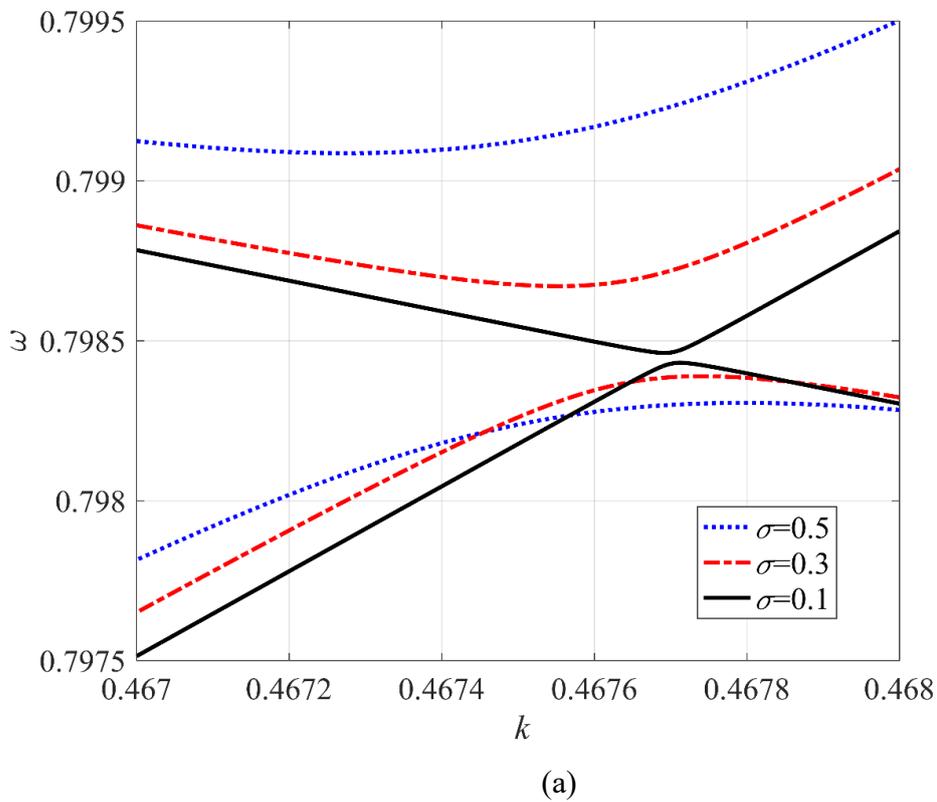

(a)

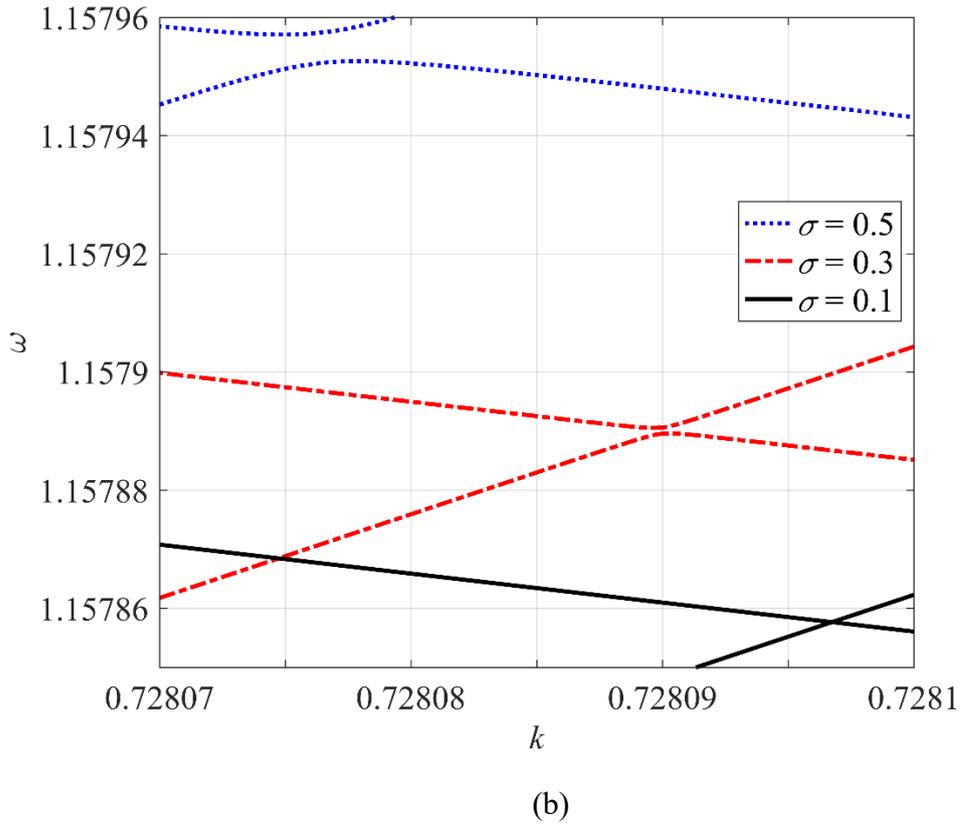

(b)

**Fig. 5** The partial enlarged view of the third band gap (a) and the fourth band gap (b) of the axially moving belt with parameters $s=0.1$, $v=0.5$ and the stiffness modulation $\sigma = 0.1$, $\sigma = 0.3$ and $\sigma = 0.5$, respectively.

## 5. Tuning band gaps by varying parameters

The dependence of the band gaps on the system parameters has been discussed in the preceding sections. If the system parameters vary, the band gaps change correspondingly, which leads to possible tunable band gap. The band gap can be manipulated by varying the parameters, such as the stiffness modulation and the foundation stiffness, which make the vibration frequency fall in the band gaps. Considering the influence of various parameters on the band gaps, the relations of the band gaps with the foundation stiffness and the stiffness modulation will be investigated in this section.

### 5.1 Varying the foundation stiffness

Considering the influence of the foundation stiffness on the band gaps, by controlling the frequency in the interval of the band gap, the wave propagation is stopped and the vibration is suppressed. For the stiffness modulation is small ($\sigma < 1$), the influence is neglectable in the first band gap. To calculate the first band gap in the first Brillouin zone, it is reasonable to take into account only one term of $f$ ($m=0$), when s<1, $\sigma < 1$. As discussed in **Section 4.1**, according to Eqs. (15) and (16), when the foundation stiffness vanishes, the group velocities are $c_g = v \pm 1$, as shown in **Fig. 2** (a). The wave propagations show in two opposite directions. When the group velocity $d\omega/dk$ vanishes, i.e., the group velocity is zero at the cut-off frequencies, then the wave propagation is stopped [29]. Thus, the wave attenuation occurs at the veering point

$$k = v\sqrt{\frac{s}{1-v^2}} . \qquad (17)$$

The critical frequency with respect to stiffness and axial velocity is obtained explicitly

$$\omega_c = \sqrt{s(1-v^2)} , \qquad (18)$$

the first band gap is below the critical frequency $\omega_c$. The width of the first band gap is presented as $|\omega_c|$. As shown in **Fig. 6**, the cut-off frequency of the first band gap decreases with the axial velocity, but increases with the foundation stiffness.

Here, the author gives an example to illustrate the vibration suppression by tuning the foundation stiffness to achieve the tunable band gap. The initial state of the system with the frequency $\omega_1$, axial velocity $v_1$ and the mean foundation stiffness $s_1$ induces a band gap. When the axial velocity varies from $v_1$ to $v_2$, the corresponding frequency will change from $\omega_1$ to $\omega_2$. If one wants to maintain the system still in the band gap for the new axial velocity $v_2$ without changing the value of the frequency, the stiffness can be tuned. As demonstrated in **Fig. 6**, the critical frequency remains stationary in the first band gap, i.e., when the axial velocity varies from $v_1$ to $v_2$, and the foundation stiffness varies from $s_1$ to $s_2$, the frequencies $\omega_1 = \omega_2$. With the aid of Eq.(15), the relation of the velocity and the stiffness can be computed as

$$\Delta s = \frac{v_2^2 - v_1^2}{1 - v_2^2} s_1 , \qquad (19)$$

where $\Delta s = s_2 - s_1$. Thus, the vibration of the system caused by the velocity perturbation can be suppressed by varying the foundation stiffness, and vice versa, if the vibration originates from the foundation stiffness, the band gap can be tuned by axial velocity

$$\Delta v = \frac{-v_1 s_2 - \sqrt{s_2^2 - (1-v_1^2)s_1 s_2}}{s_2}, \tag{20}$$

where $\Delta v = v_2 - v_1$.

Considering more than one band gap, the author will find that the foundation stiffness has significant influence on the first band gap as shown in **Fig. 7**. The cut-off frequency and the width of the band gap both increase with the foundation stiffness. The band gap will vanish if the foundation stiffness vanishes which can be drawn by comparing **Fig. 2(a)** and **Fig. 4 (a)**. According to Eq.(3), the dimensionless foundation stiffness $s$ increases with the structure period $\Phi$. If $\Phi$ decreases toward zero, the foundation stiffness decreases toward zero which means that the properties of band gap vanishes. Moreover, in contrast with the case of $\Phi$ increasing, the foundation stiffness $s$ increasing, which means that the first band gap increases with the structure period. The band gap increases with the foundation stiffness, but decreases with the axial velocity, as seen in **Fig. 8**. Both the cut-off frequency and the band gap decrease with the velocity, but the width of second band gap is not sensitive to it. Note that, there are more than two band gaps in **Fig. 8**, if one magnifies the figures, the tiny higher band gaps will emerge.

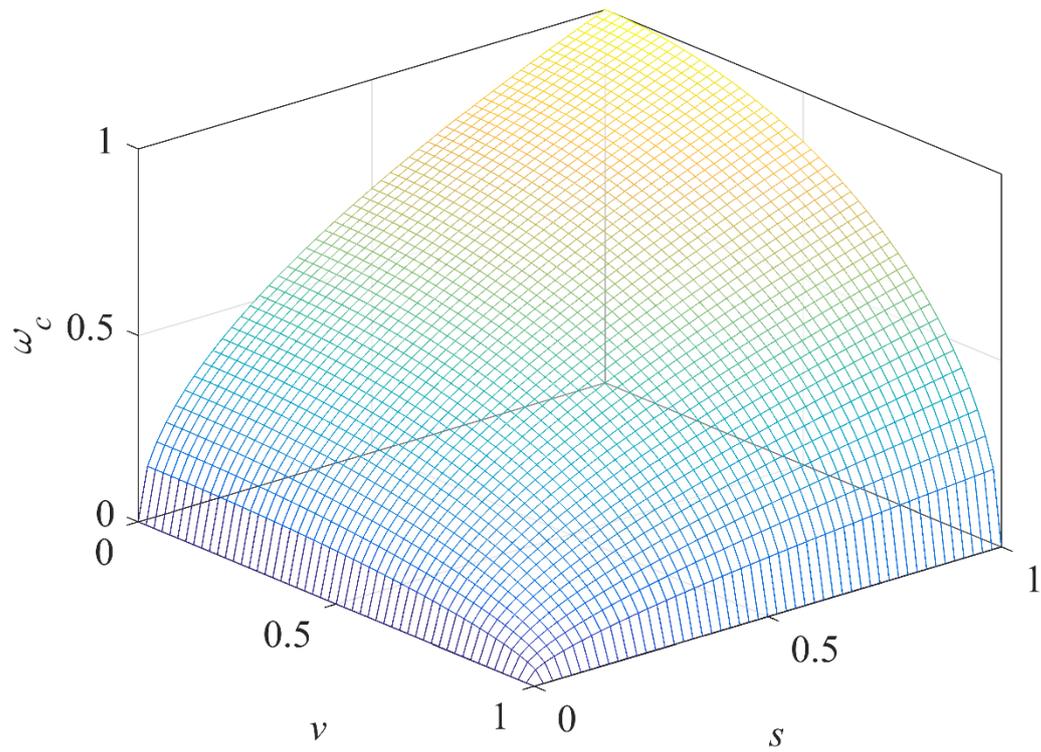

**Fig. 6** The dependent surface of critical frequency on the axial velocity and the foundation stiffness in the first band gap.

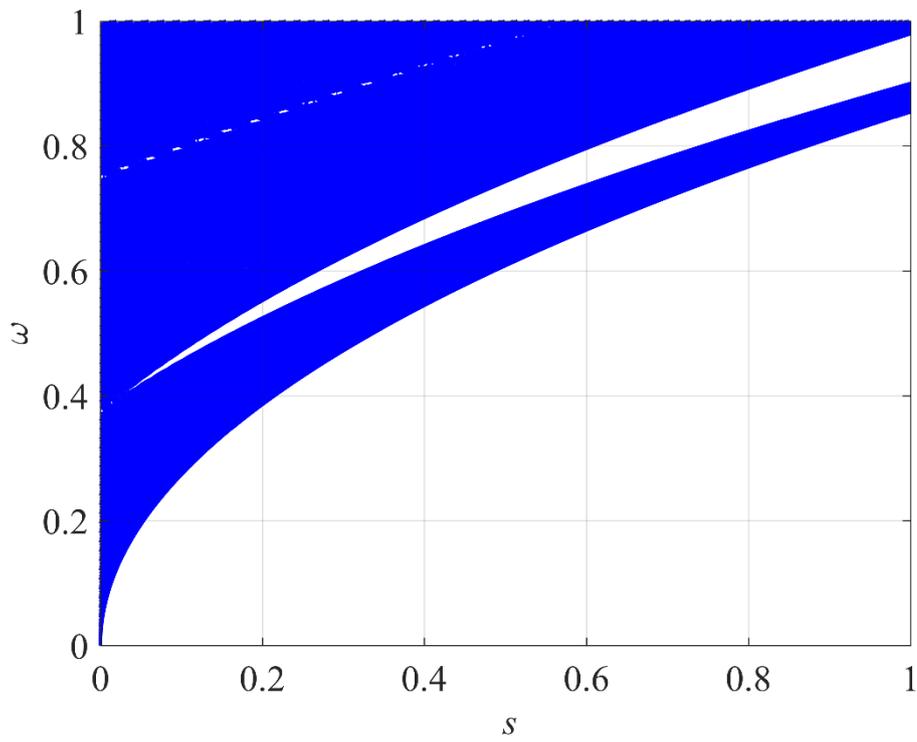

(a)

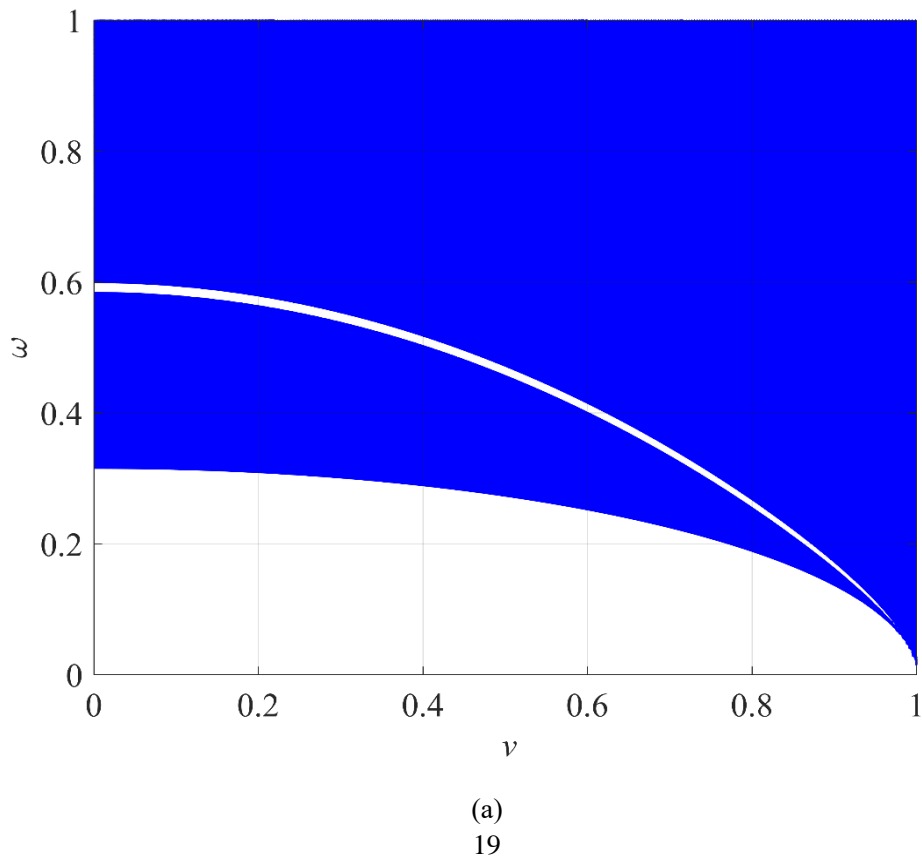

(b)

**Fig. 7** Frequency vs. foundation stiffness. The band gaps are blank areas. (a) $v$=0.5, $\sigma$=0.2; (b) $v$=0.5, $\sigma$=0.5.

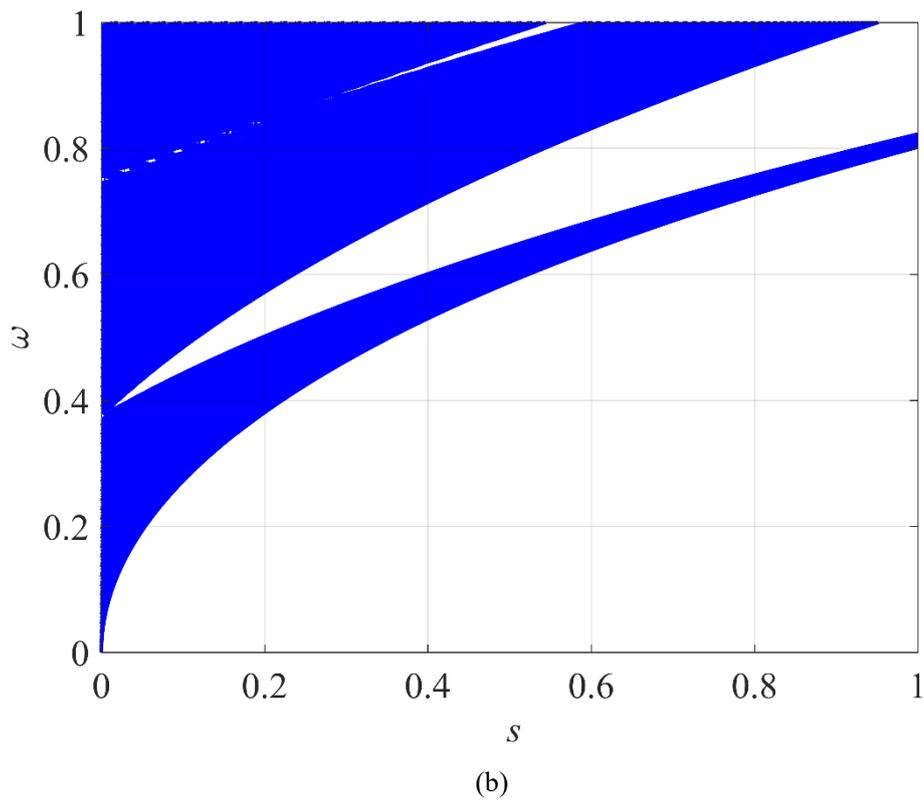

(a)

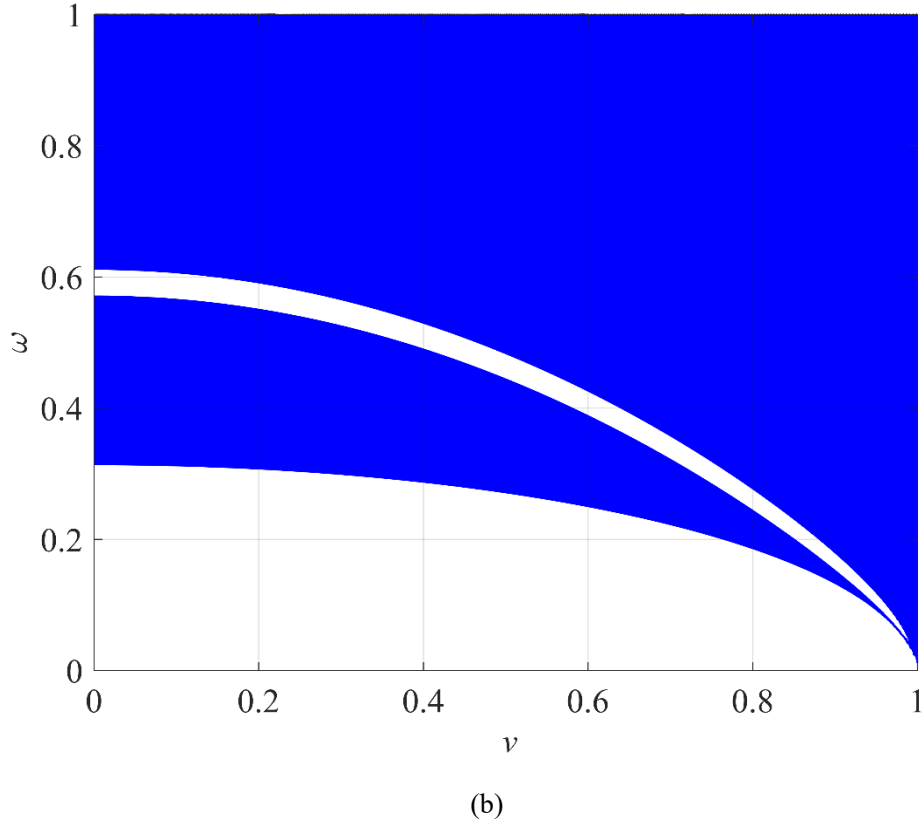

(b)

**Fig. 8** Frequency vs. axial velocity. The band gaps are blank areas. (a) $s=0.1$, $\sigma=0.2$ and (b) $s=0.1$, $\sigma=0.5$.

**5.2 Varying the stiffness modulation**

The band gap also can be manipulated by varying the stiffness modulation. To obtain the second band gap, the author considers $f$ with two terms $m=0$ and 1, and it follows that

$$f = \left((vk-\omega)^2 - k^2 - s\right)\left((v(k-1)-\omega)^2 - (k-1)^2 - s\right) - \frac{1}{4}s^2\sigma^2 = 0. \quad (21)$$

$f$ is an implicit function of $\omega$ and $k$. By applying the derivation of the implicit function, the group velocity is obtained

$$c_g = \frac{d\omega}{dk} = -\left(\frac{df}{dk}\right)\Big/\left(\frac{df}{d\omega}\right). \quad (22)$$

Similar to the previous subsection, the wave attenuation occurs at the veering point ($\omega$, $k$), when $c_g = 0$. The band gaps are obtained numerically as depicted in **Fig. 7**, **Fig. 8** and

**Fig. 9**.

The second band gap is more sensitive to the stiffness than the first one. The lower cut-off frequency of the second band gap decreases with the stiffness modulation while the upper cut-off frequency increases with it. Thus, the width of the second band gap increases with the stiffness modulation. Assuming the stiffness modulation $\sigma<1$, the first band gap almost remains stationary when the stiffness modulation changes, as shown in **Fig. 9**. In contrast with the case of band gaps illustrated in **Fig. 7**, which demonstrates that the first band gap depends more on the foundation stiffness, while the width of the second one is almost independent of the axial velocity, as depicted in **Fig. 8**. Note that, with the increasing of the stiffness modulation, the third band gap emerges, as shown in **Fig. 9** (b).

Based on the analysis above, if the frequency need to be suppressed locates in the range of the second band gap, the effect of tuning the stiffness modulation is much better than tuning the foundation stiffness. Hence, by fine design of the periodic foundation, the vibration of the axially moving belt can be suppressed for specific frequency.

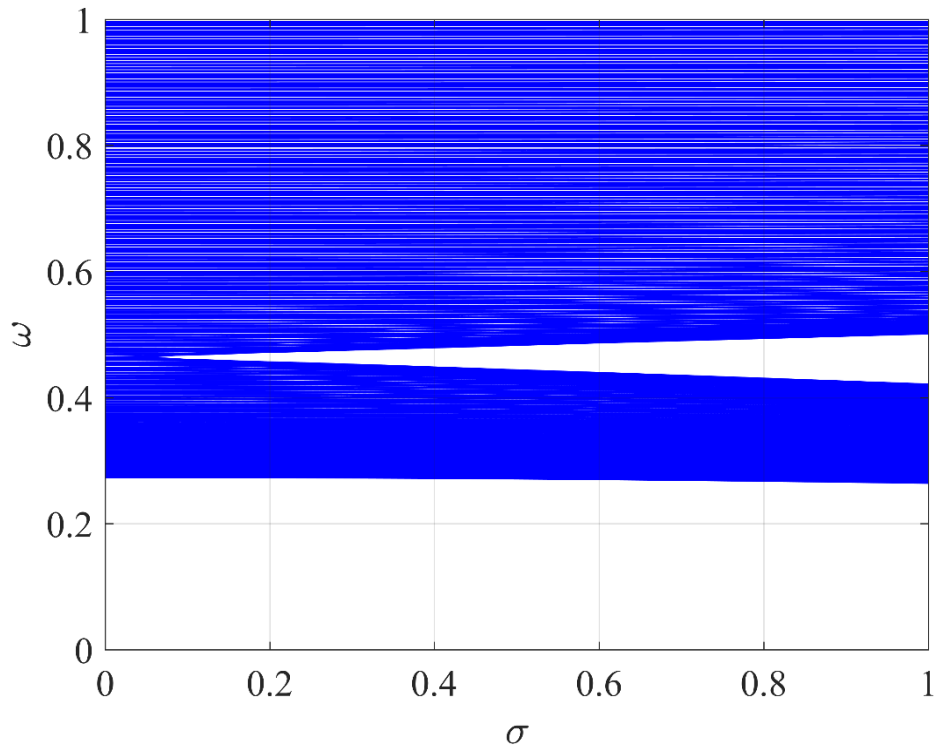

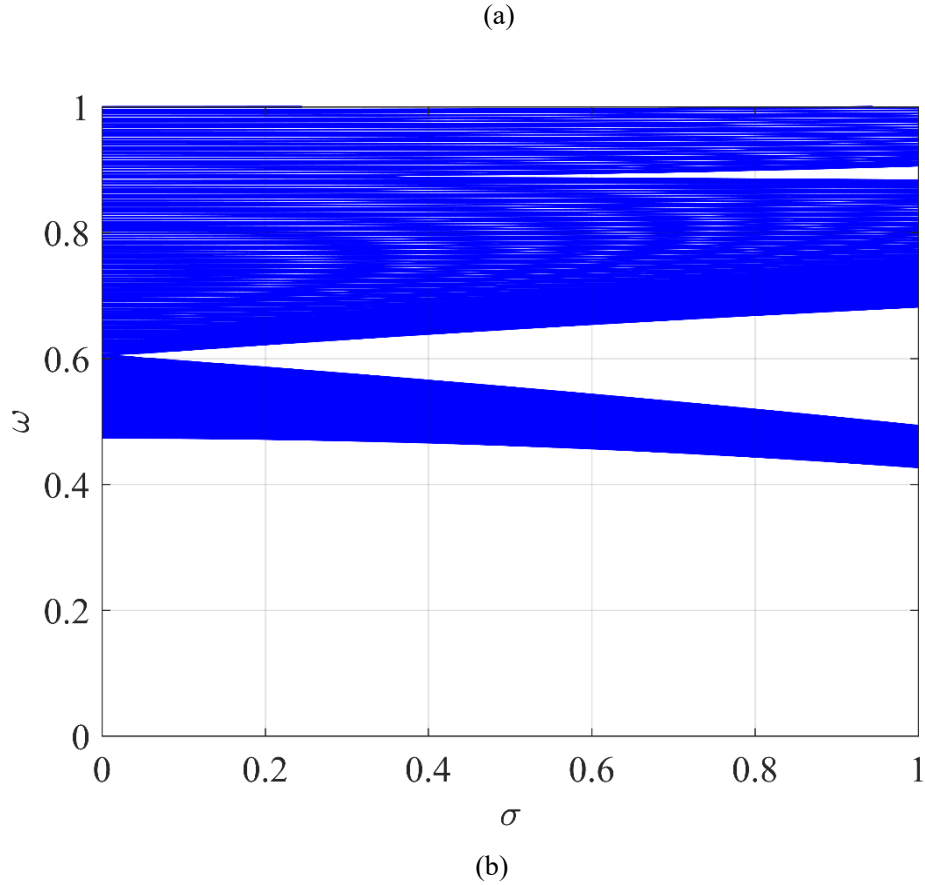

(a)

(b)

**Fig. 9** Frequency vs. stiffness modulation. The band gaps are blank areas. (a) $v$=0.5, $s$=0.1 and (b) $v$=0.5, $s$=0.3.

By making a comprehensive comparison of **Fig. 7**, **Fig. 8** and **Fig. 9**, the foundation stiffness dominates the first band gap and the stiffness modulation dominates the second and higher band gaps. Both the band gaps and the cut-off frequency increase with the foundation stiffness but decrease with the axial velocity. It should be noted that the lower cut-off frequency of the second band gap decreases with the stiffness modulation and the upper cut-off frequency increases with the stiffness modulation, which leads to the width of the second band gap increasing with the stiffness modulation.

## 6. Conclusions

This paper is devoted to describing the mechanism of the tunable band gap of the axially moving belt resting on the periodic elastic foundation. By applying the method of varying amplitudes to the periodic structure system, the relation of the band gap with the system parameters, such as stiffness modulation, foundation stiffness and the axial velocity, are obtained. Based on the relation of the parameters, the band gaps can be manipulated actively by tuning the axial velocity, the foundation stiffness or the stiffness modulation. Thus, the programmable structure can be designed according to relations of the system parameters. Three main conclusions are as follows.

(1) The wave propagation is stopped in the band gaps because the group velocity of the travelling wave is equal to zero at the veer points where the wave attenuation occurs. The eigenvalues divergence and the mode shapes veering at veer points are two different manifestations of the wave attenuation of the system.

(2) The presence of band gaps depends on the foundation stiffness. The band gaps vanish when the foundation stiffness approaches zero, but the higher band gaps more than second further depends on the stiffness modulation. The foundation stiffness dominates the first band gap, and the stiffness modulation dominates the higher band gaps.

(3) The width of band gaps increases with the foundation stiffness and the stiffness modulation, but decreases with the axial velocity. The cut-off frequency of the first band gap increases with foundation stiffness. However, the lower cut-off frequency of the second band gap decreases with the stiffness modulation, while the upper one increases with the stiffness modulation.

Based on the conclusions of band gaps, any frequency vibrations of the axially moving belts can be suppressed by designed foundation with varying stiffness. The current study may trigger new ideas to lower the vibrations of moving structures.

**Acknowledgements**

The author gratefully acknowledges the support from the National Natural Science Foundation of China (Grants No. 11972050, 11832002) and Beijing Natural Science Foundation (Grant No. 3172003).

**References**

[1] John Aitken, F. R. S. E., 1878, "An account of some experiments on rigidity produced by centrifugal force," Lond.Edinb.Phil.Mag., 5(29), pp. 81-105.
[2] Wickert, J. A., and C. D. Mote, J., 1990, "Classical Vibration Analysis of Axially Moving Continua," ASME J. Appl. Mech., 57(3), pp. 738-744.
[3] Ghayesh, M. H., and Balar, S., 2010, "Non-linear parametric vibration and stability analysis for two dynamic models of axially moving Timoshenko beams," Appl. Math. Model., 34(10), pp. 2850-2859.
[4] Yang, X. D., and Zhang, W., 2014, "Nonlinear dynamics of axially moving beam with coupled longitudinal–transversal vibrations," Nonlinear Dyn., 78(4), pp. 2547-2556.
[5] Ding, H., and Chen, L. Q., 2019, "Nonlinear vibration of a slightly curved beam with quasi-zero-stiffness isolators," Nonlinear Dyn., 95(3), pp. 2367-2382.
[6] Chen, L. Q., 2005, "Analysis and control of transverse vibrations of axially moving strings," Appl. Mech. Rev., 58(2), pp. 91-116.
[7] Ghayesh, M. H., 2009, "Stability characteristics of an axially accelerating string supported by an elastic foundation," Mech. Mach. Theory, 44(10), pp. 1964-1979.
[8] Ding, H., Li, Y., and Chen, L. Q., 2018, "Nonlinear vibration of a beam with asymmetric elastic supports," Nonlinear Dyn.(3), pp. 1-12.
[9] Ding, H., Lim, C. W., and Chen, L. Q., 2018, "Nonlinear vibration of a traveling belt with non-homogeneous boundaries," J. Sound Vib., 424(2018), pp. 78-93.
[10] Yang, X. D., Lim, C. W., and Liew, K. M., 2010, "Vibration and stability of an axially moving beam on elastic foundation," Adv. Struct. Eng., 13(2), pp. 241-248.
[11] Bhat, R. B., Xistris, G. D., and Sankar, T. S., 1982, "Dynamic behavior of a moving belt supported on elastic foundation," J. Mech. Des., 104(1), p. 143.
[12] Sorokin, V. S., and Thomsen, J. J., 2015, "Vibration suppression for strings with distributed loading using spatial cross-section modulation," J. Sound Vib., 335, pp. 66-77.
[13] Martínez-Sala, R., Sancho, J., Sánchez, J. V., Gómez, V., Llinares, J., and Meseguer, F., 1995, "Sound attenuation by sculpture," Nature, 378, p. 241.
[14] Brillouin, L., 1953, Wave propagation in periodic structures, McGraw-Hill Book Company, Inc., New York.
[15] Mead, D. J., 1970, "Free wave propagation in periodically supported, infinite beams," J. Sound Vib., 11(2), pp. 181-197.
[16] Asfar, O. R., and Nayfeh, A. H., 1983, "The application of the method of multiple scales to wave propagation in periodic structures," SIAM Rev., 25(4), pp. 455-480.
[17] Sorokin, V. S., and Thomsen, J. J., 2017, "Wave propagation in axially moving periodic strings," J. Sound Vib., 393, pp. 133-144.
[18] Sorokin, V. S., 2019, "On the effects of damping on the dynamics of axially moving spatially periodic strings," Wave Motion, 85, pp. 165-175.
[19] Yablonovitch, E., 1987, "Inhibited spontaneous emission in solid-state physics and electronics," Phys. Rev. Lett., 58(20), pp. 2059-2062.
[20] Shkunov, M. N., Vardeny, Z. V., DeLong, M. C., Polson, R. C., Zakhidov, A. A., and Baughman, R. H., 2002, "Tunable, gap-state lasing in switcbable directions for opal photonic crystals," Adv. Funct. Mater., 12(1), pp. 21-26.
[21] Fan, Q., Wang, H., Zhang, W., Wei, M., Song, Y., Zhang, W., and Yun, S., 2019, "Si–Ge alloys in C2/c phase with tunable direct band gaps: A comprehensive study," Curr. Appl. Phys., 19(12), pp. 1325-1333.

[22] Li, Z. N., Wang, Y. Z., and Wang, Y. S., 2020, "Tunable nonreciprocal transmission in nonlinear elastic wave metamaterial by initial stresses," Int. J. Solids Struct., 182-183, pp. 218-235.

[23] Airoldi, L., and Ruzzene, M., 2011, "Design of tunable acoustic metamaterials through periodic arrays of resonant shunted piezos," New J. Phys., 13, p. 21.

[24] Harrison, T. R., Hornig, G. J., Huang, C., Bu, L., Haluza-Delay, T., Scheuer, K., and Decorby, R. G., 2019, "Widely tunable bandpass filter based on resonant optical tunneling," Opt. Express, 27(16), pp. 23633-23644.

[25] Zhou, W. J., Muhammad, Chen, W. Q., Chen, Z. Y., and Lim, C. W., 2019, "Actively controllable flexural wave band gaps in beam-type acoustic metamaterials with shunted piezoelectric patches," Eur. J. Mech. a-Solid., 77, p. 11.

[26] Grinberg, I., Vakakis, A. F., and Gendelman, O. V., 2018, "Acoustic diode: Wave non-reciprocity in nonlinearly coupled waveguides," Wave Motion, 83, pp. 49-66.

[27] Lv, X. F., Chuang, K. C., and Erturk, A., 2019, "Tunable elastic metamaterials using rotatable coupled dual-beam resonators," J. Appl. Phys., 126(3), p. 9.

[28] Yang, X. D., Cui, Q. D., and Zhang, W., 2020, "Wave manipulation of two-dimensional periodic lattice by parametric excitation," ASME J. Appl. Mech., 87(1), p. 011008.

[29] Achenbach, J. D., 1975, Wave propagation in elastic solids, North-Holland Publishing Company, Amsterdam.

[30] Pierre, C., 1988, "Mode localization and eigenvalue loci veering phenomena in disordered structures," J. Sound Vib., 126(3), pp. 485-502.

[31] Leissa, A. W., 1974, "On a curve veering aberration," Journal of Applied Mathematics and Physics (ZAMP), 25(1), pp. 99-111.